\documentclass[a4paper, 11pt]{amsart}

\usepackage{amsmath,amssymb,amsthm,amscd,enumerate,setspace}
\usepackage[all]{xy}
\usepackage[dvips]{graphicx}

\usepackage{hyperref}

\input xypic
\xyoption{all}

\def\A{{\mathcal A}}
\def\a{{\bf \mathfrak a}}
\def\ann{\mbox{\rm ann }}
\def\B{{\mathcal B}}
\def\bred{\mbox{\rm bred}}
\def\br{\mbox{\rm br}}
\def\codim{\mbox{\rm codim }}
\def\coker{\mbox{\rm coker}}
\def\depth{\mbox{\rm depth }}
\def\Deg{\mbox{\rm Deg}}
\def\ddet0{\mbox{\rm det}_0 }
\def\ds{\displaystyle}
\def\E{{\mathcal E}}
\def\Ext{\mbox{\rm Ext}}
\def\Fitt{\mbox{\rm Fitt}}
\def\G{{\mathcal G}}
\def\gr{\mbox{\rm gr}}
\def\grade{\mbox{\rm grade }}
\def\Hom{\mbox{\rm Hom}}
\def\height{\mbox{\rm ht}}
\def\image{\mbox{\rm Im}}
\def\kernel{\mbox{\rm ker}}
\def\lar{\longrightarrow}
\def\m{{\mathfrak m}}
\def\Min{\mbox{\rm Min}}
\def\mult{\mbox{\rm mult}}
\def\O{{\mathcal O}}
\def\Quot{\mbox{\rm Quot}}
\def\R{{\mathcal R}}
\def\rank{\mbox{\rm rank}}
\def\rar{\rightarrow}
\def\red{{\rm  red}}
\def\reg{\mbox{\rm reg}}
\def\S{{\mathcal S}}
\def\spec{\mbox{\rm spec }}
\def\Spec{\mbox{\rm Spec}}
\def\Supp{\mbox{\rm Supp}}
\def\Sym{\mbox{\rm Sym}}
\def\T{{\mathcal T}}
\def\Tor{{\rm Tor}}
\def\tr{\mbox{\rm tr}}
\def\tratto{\mbox{\rule{2mm}{.2mm}$\;\!$}}
\def\U{{\mathcal U}}

\newtheorem{Theorem}{Theorem}[section]
\newtheorem{Lemma}[Theorem]{Lemma}
\newtheorem{Corollary}[Theorem]{Corollary}
\newtheorem{Proposition}[Theorem]{Proposition}

\newtheorem{Notation}[Theorem]{Notation}

\begin{document}

\title{Specialization and Integral Closure}

\thanks{AMS 2010 {\em Mathematics Subject Classification}. 
Primary 13B22; Secondary 13B21, 13A30. 
The first author is partially supported by the Sabbatical Leave Program at Southern Connecticut State University (Spring 2014) and the second author is partially supported by the NSF}

\author{Jooyoun Hong}
\address{ Department of Mathematics \\
Southern Connecticut State University \\
501 Crescent Street, New Haven, CT 06515}
\email{hongj2@southernct.edu}

\author{Bernd Ulrich}
\address{Department of Mathematics \\
Purdue University \\ 
150 N. University Street, West Lafayette, IN 47907}
\email{ulrich@math.purdue.edu}


\begin{abstract}
We prove that the integral closedness of any ideal of height at least two is compatible with specialization by a generic element.  
This opens the possibility for proofs using induction on the height of an ideal.
Also, with additional assumptions, we show that an element is integral over a module if it is so 
modulo a generic element of the module. This turns questions about integral closures of modules into problems 
about integral closures of ideals, by means of a construction known as Bourbaki ideal.
\end{abstract}

\maketitle

In this paper we prove in a rather general setting that the integral
closure of ideals and modules is preserved under specialization
modulo generic elements. Recall that the {\it integral closure}
$\overline I$ of an ideal $I$ in a commutative ring $R$ is the set
of all elements $y$ that are $\it integral$ over $I$, i.e., satisfy
a polynomial equation of the form
\[y^m+\alpha _1y^{m-1}+\cdots+\alpha_iy^{m-i}+\cdots +\alpha_m=0 \, , \]
where $\alpha_i \in I^i$. Alternatively, one can consider the {\em
Rees algebra} $\R (I)$ of $I$, which is the subalgebra $R[It]$ of
the polynomial ring $R[t]$. The integral closure $\overline{\R(I)}$
of $\R(I)$ in $R[t]$ is again a graded algebra, and its graded
components recover the integral closures of all powers of $I$,
\[ \overline{\R(I)} = R \oplus \overline{I}t \oplus \overline{I^2}t^2 \oplus
\  \cdots \  \oplus \overline{I^i}t^i \oplus \  \cdots \  .\] In our main result,
Theorem~\ref{extension}, we consider a Noetherian ring $R$ such that the
completion ${(R_{\m}/{\sqrt{0}})}\widehat{\;}\, $ is reduced and
equidimensional for every maximal ideal $\m$ of $R$; for instance, $R$ could be an equidimensional excellent local ring. For any
$R$--ideal $I=(a_1, \ldots, a_n)$ of height at least 2, we prove that
\[ \overline{I'/(x)}=\overline{I'}/(x) \, , \]
where $x=\sum_{i=1}^{n} z_ia_i$ is a generic element for $I$ defined
over the polynomial ring $R'=R[z_1, \ldots, z_n]$ and $I'$ denotes
the extension of $I$ to $R'$. This result can be paraphrased by
saying that an element is integral over $I$ if it is so modulo
a generic element of the ideal. Other, essentially unrelated,
results about lifting integral dependence have been proved by
Teissier, Gaffney, and Gaffney and Kleiman (\cite{G}, \cite{GK}, \cite{T1}, \cite{T2}). They go by the name `Principle of
Specialization of Integral Dependence' and play an important role in
equisingularity theory.

Our main result, Theorem 2.1, opens the possibility for proofs using induction
on the height. This yields, for instance, a quick proof of Huneke's
and Itoh's celebrated result on integral closures of powers of
complete intersections (\cite{H1}, \cite{Itoh}): If
$R$ is a ring as above and $I$ an $R$--ideal generated by a regular
sequence, then $\overline{I^{n+1}}\bigcap I^n=\overline{I}I^n$ for
every $n\geq 0$ (see Theorem~\ref{Itoh}).

 The proof of Theorem 2.1 is based on a vanishing theorem for
local cohomology modules; this is natural as the obstruction to
the specialization of integral closure lies in a cohomology module.
Thus in Theorem~\ref{H2neg} we consider a ring $R$ as above, a proper $R$--ideal
$I$ of positive height, the extended Rees algebra $\A=R[It,t^{-1}]$, its integral
closure $\overline{\A}$ in $R[t,t^{-1}]$, and an $\A$--ideal $J$ of height at least 3 generated by $t^{-1}$ and homogeneous elements of positive degree;
in this setup we show that the second local cohomology module
$H^2_J(\overline{\A})$ vanishes in non-positive degrees. Theorem~\ref{extension}
about the specialization of $\overline{I}$ follows from the
vanishing of this module in degree zero; indeed,
the algebra $\overline{\A}$
recovers the integral closures of the powers of $I$ as it coincides with
$\overline{\R(I)}$ in non-negative degrees. The main
idea in the proof of Theorem~\ref{H2neg} is to pass to the graded ring associated to the
filtration of fractional powers $\{\left(\overline{I^n} \right) ^{\frac{1}{e}} | n \geq 0 \}$
and to observe that this ring is reduced for suitable $e$ (Lemma~\ref{S}(b)). One can then use the fact that the first local
cohomology module of a non-negatively graded reduced ring vanishes
in negative degrees, which translates into the vanishing of the second local
cohomology of $\overline{\A}$ in non-positive degrees.
Theorem~\ref{H2neg} and its proof essentially go back to Itoh
(\cite[Theorem 2 and Lemma 5]{Itoh}), but special care must be exercised due to
the fact that the various integral closures appearing in the proof
need not be Noetherian. Itoh uses this result to prove his theorem
about integral closures of powers of complete intersections
(\cite[Theorem 1]{Itoh}), described above as Theorem~\ref{Itoh}; his paper
\cite{Itoh} is devoted to that proof. Based on this theorem, he
gives a proof of specialization of integral closures in
\cite[Theorem 1]{Itoh1}, but only for parameter ideals in
analytically unramified local Cohen-Macaulay rings; thus our main
theorem, Theorem~\ref{extension}, was known to Itoh for the class of complete intersection ideals.

Another set of applications of Theorem 2.1 concerns integral
closures of modules. Let $R$ be a Noetherian ring and $E$ a
finitely generated $R$--module having a {\it rank}, say $e$, by which we mean
that $K \otimes_R E$ is a free module of rank $e$ over the total
ring of quotients $K$ of $R$. In this case the {\it Rees algebra}
$\R(E)$ is defined as the symmetric algebra $\S(E)$ modulo
$R$--torsion. The Rees algebra is a standard graded $R$--algebra whose
$n^{{\rm th}}$ graded component we denote by $E^n$. If $E$ is
torsion-free, one can embed $E$ into a free module $F$ with basis
$\{t_1, \ldots, t_m\}$, and we obtain $\R(E)$ as the image of the
natural map $\S(E) \rightarrow \S(F)$, or in other words, as the
$R$-subalgebra $R[E]$ of the polynomial ring $\S(F)=R[t_1, \ldots,
t_m]$. In particular, for $E=I \subset R$ an ideal, this notion of
Rees algebra coincides with the one given at the beginning of the
introduction. If on the other hand $E$ fails to have a rank, then
the definition needs to be modified considerably, see \cite{EHU}.

Once the notion of Rees algebra is in place, one can easily
introduce the concepts of reduction, integral dependence, and
integral closure of modules. Given
two submodules $U \subset V$ of $E$, we say that $V$ is {\it
integral} over $U$ {\it in} $E$ if the inclusion of subalgebras of
$\R(E)$,
\[ R[U] \, \subset  \, R[V] \]
is an integral ring extension. The module $E$ is {\it integral} over
$U$, or $U$ is a {\it reduction} of $E$, if $E$ is integral over $U$
in $E$, or equivalently, $E^{n+1}=UE^n$ for some $n=n_0$ and hence every
$n \geq n_0$. Finally,
the {\it integral closure} $ \, \overline{U}^E$ of $U$ {\it in} $E$ is
the unique largest submodule of $E$ that contains $U$ and is
integral over $U$ in $E$; it can also be described as the degree one
component of the integral closure of the subring $R[U]$ in $\R(E)$.
If $\overline{U}^E=U$ we say that $U$ is {\it integrally closed in}
$E$. Now let $R$ be a Noetherian normal ring and $U$ a finitely generated
torsion-free $R$--module having a rank $e$. We can embed $U$ into a
free module $F=R^e$ of the same rank $e$. As $\R(U)\hookrightarrow
\R(F)$ is a birational extension and $\R(F)=\S(F)$ is a
normal ring, we see that the integral closure $\overline{\R(U)}$ of
$\R(U)$ in $\R(F)$ coincides with the integral closure of $\R(U)$ in
its total ring of quotients, and hence does not depend on the
embedding $U\hookrightarrow F$ of $U$. Thus we call $\overline{U}^F$
the {\it integral closure} of $U$ and write
$\overline{U}=\overline{U}^F$. We say that $U$ is {\it integrally
closed} provided $\overline{U}=U$. As it turns out, $\overline{U^n}$
coincides with the degree $n$ component of $\overline{\R(U)}$ for
every $n$ (\cite[(A), p.434]{Rees1}), and thus $U$ is said to be {\it
normal} if $U^n$ is integrally closed for every $n$, or
equivalently, if $\R(U)$ is a normal ring.

To state our main theorem about specialization of integral closures
of modules, Theorem~\ref{HImod2}, let $R$ be a Noetherian normal ring such
that for every maximal ideal $\m$ of $R$ the completion
$\widehat{R_{\m}}$ is reduced and
equidimensional, and let $E=\sum_{i=1}^nRa_i$ be a torsion-free
$R$--module of rank $e\geq 2$ such that for every prime ideal
$\mathfrak p$ of $R$ with $\depth R_{\mathfrak{p}} \leq 2$ either
$E_{\mathfrak p}$ has a nontrivial free direct summand or $R_{\mathfrak p}$
is regular. Under these hypotheses we prove that
\[ \overline{E'/R'x}=\overline{E'}/R'x \, , \]
where $x=\sum_{i=1}^{n} z_ia_i$ is a generic element for $E$ defined
over the polynomial ring $R'=R[z_1, \ldots, z_n]$ and $E'=R'\otimes
_RE$. The significance of this theorem is that it turns questions
about integral closures of modules into problems about integral
closures of ideals, by means of a construction known as Bourbaki
ideal. In fact, factoring out the submodule generated by $e-1$
generic elements $x_1, \ldots, x_{e-1}$ for $E$, one obtains an ideal
$I=I(E)$ called the {\it generic Bourbaki ideal} of $E$. 
Here one assumes $(R, {\mathfrak m})$ is local and one works in the
localized polynomial ring $R''=R[\{z_{ij} \}]_{{\mathfrak m}R[ \{ z_{ij} \}]}$
and the module $E''=R'' \otimes_R E$. As an
immediate consequence of Theorem~\ref{HImod2} we deduce that
\[ \overline{I} = \overline{E''}/ \sum_{i=1}^{e-1} R''x_i \,  ,  \]
provided $R$ is a ring as in the theorem and $E$ is a finitely generated
torsion-free $R$--module of rank $e>0$ such that for every prime
ideal $\mathfrak p$ of $R$ with $\depth R_{\mathfrak{p}} \leq 2$
either $E_{\mathfrak p}$ has a free direct summand of rank $e-1$ or
$R_{\mathfrak p}$ is regular, see Corollary~\ref{HIBI}.

Using this approach, one immediately deduces the results of
Kodiyalam and of Katz and Kodiyalam (\cite{KK}, \cite{K}) about
integrally closed modules over two-dimensional regular local rings
from the corresponding classical theory for ideals, going back to
Zariski (\cite{H}, \cite{ZS}), see Corollary~\ref{application1}. The proof of Theorem~\ref{HImod2} requires an
adaptation of Theorem~\ref{extension} to modules, see Theorem~\ref{HImod1}, as well as a
separate treatment of the case of two-dimensional regular local
rings, see Theorem~\ref{2dimRLR}. To show the latter result we use aspects of
Kodiyalam's work, which we summarize in Theorem~\ref{transform}. For the most
part however, we recover his results through the method of Bourbaki
ideals.

\bigskip

\section{A vanishing theorem for local cohomology modules}

\smallskip

In this section we present the vanishing theorem that will be used
in the proof of our main result about specialization of integral closure.

Let $R$ be a Noetherian ring and $I$ an $R$--ideal. Denote the
extended Rees ring of $I$ by $\A=R[It,t^{-1}]$ and the integral
closure of $\A$ in $R[t,t^{-1}]$ by $\overline{\A}$. 
Fix an integer $e>0$. Let $u$ be a variable with $u^e=t$ and $\deg u = \frac{1}{e}$. Thus the $R$--algebra
$R[u,u^{-1}]$ is graded by ${\frac{1}{e}} \, {\mathbb Z}$ and it contains $R[t,t^{-1}]$
as a Veronese subring. Consider
the finitely generated graded $R$--subalgebra $\S=\A[u^{-1}]$ of $R[u,u^{-1}]$. This ring is of the form
\[ \S= \ \cdots \  \oplus Ru^{-1} \oplus R
\oplus Iu \oplus Iu^2 \oplus \ \cdots \ \oplus Iu^e \oplus I^2 u^{e+1}
\oplus I^2 u^{e+2} \oplus \  \cdots \  \oplus I^2u^{2e} \oplus \  \cdots \ .
\]
Notice that $\A$ is a Veronese subring of $\S$.
Denote the integral closure of $\S$ in $R[u,u^{-1}]$ by
$\overline{\S}$. The ring $\overline{\S}$ is a graded $R$--subalgebra 
of $R[u, u^{-1}]$ and $\overline{\A}$ is a  Veronese subring of
$\overline{\S}$ (\cite[2.3.2]{HSw}). Note that $\overline{\S}$ and $\overline{\A}$ 
need not be Noetherian though. Let $a$ be an element of $R$ and $n$ an integer.
Then $au^n \in \overline{\S}$ if and only if $a^et^n \in
\overline{\A}$, or equivalently, $a^e \in \overline{I^n}$, where $I^n=R$ if $n \leq 0$.
Write $\left( \overline{I^n} \right)^{\frac{1}{e}}=\{a \in R\;|\; a^e \in \overline{I^n} \}$.
Thus 
$\overline{\S}$ is of the form
\[\overline{\S}=  \bigoplus_{n \in {\mathbb Z}} \left( \overline{I^n} \right)^{\frac{1}{e}} u^n \, .\] 

Finally, for a ring $R$ with nilradical $\sqrt{0}$, we denote
$R/\sqrt{0} \, $ by $R_{\red}$. We say that a Noetherian local ring is {\it analytically unramified}
if its completion $\widehat{R}$ is reduced and that a Noetherian ring is {\it locally analytically 
unramified} if each of its localizations at a maximal ideal is analytically unramified.  

In the next lemma we prove that for a suitable choice of $e$, the associated graded ring 
\[ \E = \overline{\S}/u^{-1}\overline{\S}= \bigoplus_{n \in {\mathbb Z}_{\geq 0}}  
\left( (\overline{I^n} )^{\frac{1}{e}} /   ( \overline{I^{n+1}} )^{\frac{1}{e}} \right) \] is reduced. From
this we deduce that $\E$ is Noetherian and that its integral closure in the total ring of 
quotients is non-negatively graded, though not necessarily Noetherian.

\begin{Lemma}\label{S}
Let $R$ be a Noetherian, equidimensional, universally catenary local ring of dimension $d$ such that  $R_{\red}$ is
analytically unramified. Let $I=(a_1, \ldots, a_n)$ be a proper $R$--ideal
with $\height(I)>0$ and write $\A=R[It,t^{-1}]$ for the extended Rees ring of
$I$. Let $\upsilon_1, \ldots, \upsilon_r$ be the Rees valuations of $I$
$($\cite[10.1.1]{HSw}$)$, and let $e$ be the least common multiple of the values
$\upsilon_1(I), \ldots, \upsilon_r(I)$ of $I$. Let $u$ be a variable with $u^e=t$ and $\deg u = {\frac{1}{e}}$. 
Write $\S=\A[u^{-1}]$ and let $\overline{\S}$ be the integral closure of
$\S$ in $R[u,u^{-1}]$.
\begin{enumerate}
\item[{\rm (a)}] Let  $\,\overline{\S_{\red}}$ denote the integral closure of $\S_{\red}$ in
$R_{\red}[u,u^{-1}]$. Then the $R$--algebra $\,\overline{\S_{\red}}$ is finitely
generated and graded by ${\frac{1}{e}} \, {\mathbb Z}$, has a unique maximal homogeneous ideal, which is a maximal ideal, and is
equidimensional of dimension $d+1$.

\item[{\rm (b)}] One has the equality $\overline{\S}/u^{-1}\overline{\S}=\overline{\S_{\red}}/u^{-1}\overline{\S_{\red}} \,$.
This $R$--algebra is finitely generated and graded by ${\frac{1}{e}} \, {\mathbb Z}_{\geq 0}$, has a unique maximal homogeneous ideal,
is equidimensional of dimension $d$, and is reduced.

\item[{\rm (c)}] Let $R'=R[z_1, \ldots, z_n]$ be a polynomial ring in the variables $z_1, \ldots, z_n$, $I'=IR'$, and
$x=\sum_{i=1}^n z_i a_i$. Then $t^{-1}, xt$ is a regular sequence on $R' \otimes_R \overline{\A}$. In particular, $xt$ is
regular on $R' \otimes_R(\overline{\A}/t^{-1}\overline{\A})$.

\item[{\rm (d)}] Write $\E=\overline{\S}/u^{-1}\overline{\S}$ and let $\overline{\E}$ be the integral closure of
$\E$ in its total ring of quotients. Then $\overline{\E}$ is a finite direct product of Krull domains that are
graded by ${\frac{1}{e}} \, {\mathbb Z}_{\geq 0} \, $. In particular, $\overline{\E}$ is
graded by ${\frac{1}{e}} \, {\mathbb Z}_{\geq 0}$.
\end{enumerate}
\end{Lemma}

\begin{proof} (a)  
Let $I_{\red}$ be the image of $I$ in $R_{\red}$. We denote the degree $m$
component by $[\;\;]_m$. Then
\[ [ \overline{\S_{\red}}]_{{\frac{n}{e}}} \cdot u^{-n} = \left( \overline{I_{\red}^{n}} \right)^{\frac{1}{e}} \subset   
\left( \overline{ I_{\red}^{ \lfloor \frac{n}{e} \rfloor e  } } \right)^{\frac{1}{e}} 
= \overline{ I_{\red}^{ \lfloor \frac{n}{e}  \rfloor}  }  \subset I_{\red}^{\lfloor \frac{n}{e} \rfloor  -k} \, ,
\] 
where the last containment holds for some fixed integer $k \geq 0$ and every integer $n$ because $R_{\red}$ is analytically unramified
(\cite[1.4]{Rees0}). On the other hand,
\[ I_{\red}^{\lfloor \frac{n}{e} \rfloor  -k} = [ \S_{\red} ]_{( \lfloor \frac{n}{e}  
\rfloor -k   ) }  \cdot u^{ - \left( \lfloor \frac{n}{e}  \rfloor -k \right) e } 
\subset 
[ \S_{\red} ]_{( {\frac{n}{e}}-(k+1) )} \cdot u^{-n+(k+1)e }.   \]
Therefore $ \, {\ds \overline{\S_{\red}} \subset u^{(k+1)e} \, \S_{\red}}$. 
Since $\S_{\red}$ is a Noetherian ring, it follows that $\overline{\S_{\red}}$ is a finitely generated
$\S_{\red}$--module, hence a finitely generated $R_{\red}$--algebra. We already argued that  
$\,\overline{\S_{\red}}$ is graded by ${\frac{1}{e}} \, {\mathbb Z}$.
Writing ${\mathfrak m}$ for the maximal ideal of $R$, we have $(\overline{I})^{\frac{1}{e}} \subset {\mathfrak m}$.
Therefore $( \overline{\S_{\red}})_{<0} \oplus {\mathfrak m} \oplus ( \overline{\S_{\red}})_{>0}$ is a maximal ideal 
of $\overline{\S_{\red}}$ and the unique maximal homogeneous ideal.

Since $\overline{\S_{\red}}$ is torsion-free
over $\A_{\red}$ and both rings are Noetherian and reduced, for any
minimal prime ideal $\mathfrak{Q}$ of
$\overline{\S_{\red}}$, the ideal
$\mathfrak{q}=\mathfrak{Q} \cap \A_{\red}$ is a minimal prime ideal of
$\A_{\red}$. The extension $\A_{\red}/\mathfrak{q} \subset \overline{\S_{\red}}/\mathfrak{Q} \, $
is integral and the ring $\A_{\red}$ is equidimensional of dimension $d+1$. Therefore
\[ \dim(\overline{\S_{\red}}/\mathfrak{Q})=\dim(\A_{\red}/\mathfrak{q})= d+1 \, , \]
from which we deduce that $\overline{\S_{\red}}$ is equidimensional of dimension $d+1$.

\medskip

\noindent (b) We first show that $u^{-1}\overline{\S}$ is a radical ideal. Thus let $au^n \in \overline{\S}$ 
be a homogeneous element such that
$(au^n)^m \in u^{-1}\overline{\S}$ for some positive integer $m$.
Then $a^{me} \in \overline{I^{nm+1}}$. Hence for each $i=1, \ldots,
r$, we have ${\ds \frac{e}{\upsilon_i(I)} \upsilon_i(a) \geq n+
\frac{1}{m}}$. Since $e/\upsilon_i(I)$ is an integer for each $i$,
we obtain ${\ds \frac{e}{\upsilon_i(I)} \upsilon_i(a) \geq n+ 1}$, which means
$a^e \in \overline{I^{n+1}}$. Therefore $au^{n+1} \in
( \overline{I^{n+1}} )^{\frac{1}{e}}   u^{n+1} \subset \overline{\S}$
and hence $au^{n} \in  u^{-1}\overline{\S}$. Thus we have proved that $u^{-1}\overline{\S}$ is a radical ideal.

Therefore $\overline{\S}/u^{-1}\overline{\S}$ is reduced. Since moreover $(\overline{\S})_{\red} = \overline{\S_{\red}}$, we conclude 
that $\overline{\S}/u^{-1}\overline{\S}=\overline{\S_{\red}}/u^{-1}\overline{\S_{\red}}$. 
The remaining assertions follow from part (a) because $u^{-1}$ is a homogeneous regular element on $\overline{\S_{\red}}$
and $\overline{\S_{\red}}$ is a catenary ring.

\medskip

\noindent (c) Write $\E={\overline{\S}}/u^{-1}{\overline{\S}}$. From part (b) we obtain 
\[\height(\E_+)=\dim(\E)-\dim(\E/{\E_+}) = d - \dim(R/I) >0 \, .\] 
It follows that $\grade \E_+ >0$ because $\E$ is reduced. We claim that $\E_+ \subset \sqrt{I^*\E}$, 
where $I^*$ denotes the image of $I$ in $\E_1$. Any homogeneous element in $\E_{+}$ is
of the form $\alpha^* \in \E_{\frac{k}{e}}$, where $\alpha^*$ is the
image of an element $\alpha \in ( \overline{I^k} )^{\small{\frac{1}{e}}}$ and $k$ is a positive integer.
Then $\alpha^e \in \overline{I^k}$. For a sufficiently large integer $l$, we obtain $\alpha^{el} \in
\overline{I^{kl}}=I\overline{I^{kl-1}}$. Notice that
$(\alpha u^k)^{el}= \alpha^{el}u^{ekl} \in \overline{\S}_{kl}$, $Iu^e \subset
\overline{\S}_1$, and $ \overline{I^{kl-1}}u^{e(kl-1)} \subset
\overline{\S}_{kl-1}$. Hence ${\alpha^*}^{el} \in I^*\E$, which
proves the claim. From the claim we obtain $\grade I^*\E>0$. 

Since the leading form
$x^*$ of $x$ in $R' \otimes_R \E$ is a generic element for $I^* \E$, it then follows that
$x^*$ is regular on $R' \otimes_R \E$. Therefore $u^{-1}, xu^e$ is a regular sequence 
on $R' \otimes_R\overline{\S}$. We conclude that $(u^{-1})^e=t^{-1}, xt$ is a
regular sequence on $R' \otimes_R \overline{\A}$ because $R' \otimes_R \overline{\A}$ is 
a direct summand of $R' \otimes_R \overline{\S}$. In particular,
$xt$ is a regular element on
$R' \otimes_R (\overline{\A}/t^{-1}\overline{\A})$.

\medskip

\noindent (d) Since the ring $\E$ is Noetherian, it has only finitely many
minimal prime ideals, say $\mathfrak{Q}_1, \ldots,
\mathfrak{Q}_s$. Notice that $\mathfrak{Q}_1, \ldots,
\mathfrak{Q}_s$ are all homogeneous.  Let
$\overline{\E/\mathfrak{Q}_i}$ denote the integral closure of
$\E/\mathfrak{Q}_i$ in its quotient field for each $i=1, \ldots,
s$. Since $\E$ is reduced, the integral closure $\overline{\E}$
equals $\overline{\E/\mathfrak{Q}_1} \times \cdots \times
\overline{\E/\mathfrak{Q}_s}$.  By \cite[4.10.5 and 2.3.6]{HSw}, the integral
closures $\overline{\E/\mathfrak{Q}_i}$ are Krull domains and are graded by ${\frac{1}{e}} \, {\mathbb Z}_{\geq 0}$. 
In particular, $\overline{\E}$ is graded by ${\frac{1}{e}} \, {\mathbb Z}_{\geq 0}$. 
\end{proof}

\medskip

Here is the vanishing theorem that will be used in the next section.

\begin{Theorem}\label{H2neg}
Let $R$ be a Noetherian, locally equidimensional, universally
catenary ring such that $R_{\red}$ is locally analytically
unramified. Let $I$ be a proper $R$--ideal with $\height(I)>0$, $\A=R[It,t^{-1}]$ the extended
Rees ring of $I$, and $\overline{\A}$ the integral closure of $\A$
in $R[t,t^{-1}]$. Let $J$ be an $\A$--ideal of height at
least $3$ generated by $t^{-1}$ and homogenous elements of positive degree. Then $\left[ H^2_{J}(\overline{\A})
\right]_n =0$ for all $n \leq 0$, where $[\;\;]_n$ denotes the degree $n$
component.
\end{Theorem}

\begin{proof} By localizing, we may assume that $R_{\red}$ is an 
analytically unramified Noetherian local  ring.  Let $\S=\A[u^{-1}]$, 
where $u$ is a variable with $u^e=t$ and $\deg u = {\frac{1}{e}}$ and $e$ 
is the least common multiple of the  values of $I$ under its Rees valuations. 
Denote the integral closure of $\S$ in $R[u,u^{-1}]$ by $\overline{\S}$.  
It suffices to show that
$\left[H^2_{J}(\overline{\S})\right]_n =0$ for any rational number $n \leq 0$
because $\overline{\A}$ is a Veronese subring
of $\overline{\S}$. Let $d=\dim R$. Write $\E=\overline{\S}/u^{-1}\overline{\S}$ 
and let $\overline{\E}$ be the integral closure
of $\E$ in its total ring of quotients. Then $\E$ is a finitely generated $R$--algebra 
and $\overline{\E}$ is graded by  ${\frac{1}{e}} \, {\mathbb Z}_{\geq 0}$ 
(Lemma~\ref{S}(b),(d)). Some of the technical difficulties in the following 
proof stem from the fact that the ring $\overline{\E}$ may be non-Noetherian.

First we show that $\height(J\mathcal{F}) \geq 2$ for any finitely
generated graded $R$--subalgebra $\mathcal{F}$ of $\overline{\E}$ containing
$\E$. Since $\dim (\A/J) \leq d-2$, we have
$\dim(\mathcal{F}/J\mathcal{F}) \leq d-2$. Notice that $J\mathcal{F}$ is a homogeneous ideal
and that $\mathcal{F}$ is a catenary ring with a unique maximal homogeneous ideal. Thus 
it remains to show that $\mathcal{F}$ is equidimensional of dimension $d$. Since $\mathcal{F}$ 
is a torsion-free $\E$--module and $\E$ is Noetherian and reduced (Lemma~\ref{S}(b)), for any 
minimal prime ideal $\mathfrak{Q}$ of $\mathcal{F}$, the
ideal $\mathfrak{q}=\mathfrak{Q} \cap \E$ is a minimal prime of $\E$. The extension 
$\E/\mathfrak{q} \subset \mathcal{F}/\mathfrak{Q} \, $ is integral and the ring $\E$
is equidimensional of dimension $d$ (Lemma~\ref{S}(b)). Hence we obtain
$\dim(\mathcal{F}/\mathfrak{Q})=\dim(\E/\mathfrak{q}) = d$, as asserted.

Next we claim that there exists an $\overline{\E}$--regular sequence of length $2$
in $J\overline{\E}$. The integral closure
$\overline{\E}$ is a direct product of finitely many Krull domains (Lemma~\ref{S}(d)). 
Recall that a principal ideal generated by a nonzero element in a Krull domain is a 
finite intersection of primary ideals of height $1$ (\cite[4.10.3]{HSw}).
Hence by Prime Avoidance it suffices to show that $\height(J\overline{\E}) \geq 2$.
Suppose that $\height(J\overline{\E}) \leq 1$. In this 
case, $J \overline{\E}$ has a minimal prime $\mathfrak{P}$ with $\height(\mathfrak{P}) = \height( J \overline{\E} ) \leq 1$.
Since $J \overline{\E}$ is homogeneous, the prime ideal $\mathfrak{P}$ is homogeneous 
(\cite[1.5.6(a),(b.i)]{BH}). We claim that there exists a homogenous element 
$y \in J \overline{\E}$ such that $\mathfrak{P}$ is a minimal prime ideal of 
$y\overline{\E}$. If $\height(\mathfrak{P})=0$, this is obvious. 
If $\height(\mathfrak{P}) = \height (J \overline{\E} ) =1$, then $J \overline{\E}$ 
is not contained in any of the finitely many minimal prime ideals of $\overline{\E}$. 
As $J \overline{\E}$ is generated by homogeneous elements of positive degree, it 
follows that $J \overline{\E}$ contains a homogeneous element $y$ not contained 
in any minimal prime of $\overline{\E}$ (\cite[1.5.10]{BH}). Hence $\mathfrak{P}$ 
is a minimal prime ideal of $y \overline{\E}$ because $\height(\mathfrak{P})=1$.
Since the ideal $J$ is finitely generated, $y\overline{\E}_{\mathfrak{P}}$
contains $J^n \overline{\E}_{\mathfrak{P}}$ for some $n$. Let $J^n=(f_1, \ldots, f_s)$, where
$f_i$ are homogeneous. For each $i$, there exists a homogeneous
element $u_i \in \overline{\E} \setminus \mathfrak{P}$ such that
$f_iu_i \in y\overline{\E}$. Write $f_iu_i = y a_i $ for some
homogeneous elements $a_i \in \overline{\E}$. Let
$\mathcal{F}=\E[u_1, \ldots, u_s, a_1, \ldots, a_s, y]$ and set
$\mathfrak{p}=\mathfrak{P} \cap \mathcal{F}$.
Then $\mathcal{F}$ is a finitely generated graded $R$--subalgebra of
$\overline{\E}$ containing $\E$. Hence we have $\height(J\mathcal{F}) \geq 2$
as shown in the preceding paragraph. On the other hand, for every $i$, $u_i \in
\mathcal{F} \setminus \mathfrak{p}$ and $y \in \mathfrak{p}$. Since
$f_iu_i \in y \mathcal{F}$, we obtain $f_i \in y \mathcal{F}_{\mathfrak{p}}$
for every $i$. Therefore $J^n\mathcal{F}_{\mathfrak{p}}$ is contained
in $y\mathcal{F}_{\mathfrak{p}}$. But $y\mathcal{F}_{\mathfrak{p}}$
is a proper principal ideal of $\mathcal{F}_{\mathfrak{p}}$ since $y \in
\mathfrak{p}$, and the ring $\mathcal{F}_{\mathfrak{p}}$ is Noetherian. Hence $\height(J\mathcal{F}) \leq 1$, which is a
contradiction.

Now we are ready to complete the proof. From the exact sequence of ${\frac{1}{e}} \, {\mathbb Z}$--graded modules  
\[\begin{CD} 
0 \rightarrow \overline{\S}({\small{\frac{1}{e}}}) \stackrel{ u^{-1}}{\longrightarrow}
\overline{\S} \rightarrow \E \rightarrow 0  \, , 
\end{CD}
\]
we obtain a long exact sequence
\[\begin{CD}
\cdots @>>> H^1_{J}(\E) @>>> H^2_{J}(\overline{\S}({\small{\frac{1}{e}}})) @>u^{-1}>>
H^2_{J}(\overline{\S}) @>>> \cdots.
\end{CD}\]
Suppose that $a \in
\left[H^2_{J}(\overline{\S})\right]_n=\left[H^2_{J}(\overline{\S}({\frac{1}{e}}))\right]_{n-{\frac{1}{e}}}
$ for some rational number $n \leq 0$ and that $a \neq 0$. Let $l$ be the smallest
non-negative integer such that $(u^{-1})^l a=0$. Such integer $l$ exists
because $u^{-1}$ belongs to $\sqrt{J}$, and clearly $l >0$ because $a \neq 0$ . Therefore $u^{-l+1}a \in
\ker(\cdot\; u^{-1}) \cap
\left[H^2_{J}(\overline{\S}({\frac{1}{e}}))\right]_{n-\frac{l}{e}}$ and $n- \frac{l}{e} <0$. Now
if we show that $[\ker(\cdot\; u^{-1})]_m =0$ for all rational numbers $m <0$,
then $0=u^{-l+1}a=(u^{-1})^{l-1} a$, which contradicts  the
minimality of $l$. Hence it suffices to prove that $\left[
H^1_{J}(\E) \right]_m =0$ for all  rational numbers $m <0$. 

Since there exists an
$\overline{\E}$--regular sequence of length at least $2$ in
$J\overline{\E}$, we get
$H^0_{J}(\overline{\E})=H^1_{J}(\overline{\E})=0$.
The exact sequence $ 0 \rightarrow \E \rightarrow \overline{\E}
\rightarrow \overline{\E}/\E \rightarrow 0$ yields
\[\begin{CD} 0=H^0_{J}(\overline{\E}) @>>> H^0_{J}(\overline{
\E}/\E) @>\cong>> H^1_{J}(\E) @>>> H^1_{J}(\overline{\E})=0 \, .
\end{CD} \] Since $\overline{\E}$ is graded by ${\frac{1}{e}} \, {\mathbb Z}_{\geq 0}$,
hence so is $H^0_{J}(\overline{\E}/\E) \subset \overline{\E}/\E$ and
we conclude that
$\left[ H^1_{J}(\E) \right]_m \cong \left[ H^0_{J}(\overline{\E}/\E)
\right]_m = 0$ for all  rational numbers $m <0$. 
\end{proof}

\medskip

\section{Specialization by generic elements of ideals}

\smallskip

We show that the integral closure of any ideal of height at least
$2$ is compatible with specialization by generic elements. When the
ideal is a complete intersection, this theorem is proved in
\cite[Theorem 1]{Itoh1}.

\begin{Theorem}\label{extension}
Let $R$ be a Noetherian, locally equidimensional, universally
catenary ring such that ${R}_{\red}$ is locally analytically
unramified. Let $I=(a_1, \ldots, a_n)$ be an $R$--ideal of height at
least $2$. Let $R'=R[z_1, \ldots, z_n]$ be a polynomial ring in the variables
$z_1, \ldots, z_n$, $I'=IR'$, and $x={
\sum_{i=1}^{n}} z_ia_i$. Then $ \overline{I'/(x)} =
\overline{I'}/(x)$.
\end{Theorem}

\begin{proof} We may assume that $R$ is a local ring of dimension $d$ and $I$ is a proper ideal.
Let $\A$ and $\B$ be the extended Rees rings of $I'$ and
$I'/(x)$, respectively. Denote $R'/(x)$ by $\mathfrak{R}$. Let $\overline{\A}$
and $\overline{\B}$ be the integral closures of $\A$ and $\B$ in
$R'[t,t^{-1}]$ and $\mathfrak{R}[t,t^{-1}]$, respectively. 
Consider the natural map
$\varphi:  \overline{\A}/xt\overline{\A} \longrightarrow
\overline{\B}$. Let $J$ be the $\A$--ideal $(t^{-1}, I't)$, which is generated by 
$t^{-1}$ and by linear forms. The ideal $J$ has height at least $3$. This holds
because the $R'$--ideal $I'$ has height at least $2$ and the ring $\A$ is catenary and equidimensional
of dimension $d+1$ with a unique maximal homogeneous ideal, which is a maximal ideal.

First we claim that $\varphi_{\mathfrak{p}}$ is an isomorphism for
all $\mathfrak{p} \notin V(J)$. Suppose $t^{-1} \notin
\mathfrak{p}$. In this case the claim follows because $\varphi_{t^{-1}}$ factors through the natural isomorphisms
\[\left(\overline{\A}/xt\overline{\A} \right)_{t^{-1}} \cong \overline{\A}_{t^{-1}}/xt\overline{\A}_{t^{-1}} \cong
R'[t,t^{-1}]/xR'[t,t^{-1}] \cong \mathfrak{R}[t,t^{-1}] \cong
\overline{\B}_{t^{-1}} \, .
\] Suppose $I't \not \subset \mathfrak{p}$. We may assume that $a_nt
\notin \mathfrak{p}$.
One has \[\overline{\A}=\overline{R'[I't,t^{-1}]} =\overline{R[z_1,
\ldots, z_n][I't,t^{-1}]} \cong \overline{R[It,t^{-1}]}[z_1,
\ldots, z_n] \, ,\]where $\overline{R[It,t^{-1}]}$ is the integral
closure of $R[It,t^{-1}]$ in $R[t,t^{-1}]$. Therefore
\[ \begin{array}{lll}\left(\overline{\A}/xt\overline{\A} \right)_{a_nt} &\cong
&\left(\overline{R[It, t^{-1}]}[z_1, \ldots, z_n]/(xt)\right)_{a_nt}
\\& & \\
&\cong & \overline{R[It, t^{-1}]}_{a_nt}[z_1, \ldots, z_{n}]/ \left({\ds
\sum_{i=1}^{n-1}\frac{a_iz_it}{a_nt}}+z_n \right) \\ & &\\&\cong
&\overline{R[It, t^{-1}]}_{a_nt}[z_1, \ldots, z_{n-1}] \, .
\end{array}\]
On the other hand, $(\overline{\B})_{a_nt} \subset \mathfrak{R}[t,t^{-1}]_{a_nt}$ and
\[
\mathfrak{R}[t,t^{-1}]_{a_nt}  \cong  \left(R'[t,t^{-1}]/(x) \right)_{a_nt}   \cong  \left(R'[t,t^{-1}]/(xt) \right)_{a_nt}
 \cong \left(R[t,t^{-1}]_{a_nt}\right)[z_1, \ldots, z_{n-1}] \, . \] Now we
have

\[\xymatrix{
\left(\overline{\A}/xt\overline{\A} \right)_{a_nt} \ar[rr]
\ar[d]^{\cong}
& &(\overline{\B})_{a_nt}\; \ar@{^{(}->}[rr] & & \mathfrak{R}[t,t^{-1}]_{a_nt} \ar[d]^{\cong} \\
\left(\overline{R[It, t^{-1}]}_{a_nt} \right)[z_1, \ldots,
z_{n-1}]\; \ar@{^{(}->}[rrrr] &&&&
\left(R[t,t^{-1}]_{a_nt}\right)[z_1, \ldots, z_{n-1}] \, . }
\]
Since $\left(\overline{R[It, t^{-1}]}_{a_nt} \right)[z_1, \ldots,
z_{n-1}]$ is integrally closed in
$\left(R[t,t^{-1}]_{a_nt}\right)[z_1, \ldots, z_{n-1}]$, this proves
the claim.

Denote the kernel and cokernel of $\varphi$ by $K$ and $C$,
respectively. Our goal is to prove that $C_1=0$, which we do by
identifying $C$ with a submodule of $H^2_J(\overline{\A})(-1)$. By the claim above, we
have $K=H^0_J(K)$ and $C=H^0_J(C)$. Moreover, $H^0_J(\overline{\B})=0$ 
because $J$ contains the $\overline{\B}$--regular element $t^{-1}$.
From the exact sequence 
$0 \rightarrow K \rightarrow \overline{\A}/xt\overline{\A}
\stackrel{\varphi}{\longrightarrow} \image(\varphi) \rightarrow 0$,
we obtain $H^i_J(\overline{\A}/xt\overline{\A}) \cong
H^i_J(\image(\varphi))$ for all $i \geq 1$. Then,
using the exact sequence $0 \rightarrow \image(\varphi) \rightarrow
\overline{\B} \rightarrow C \rightarrow 0$, we get
\[\begin{CD} 0=
H^0_J(\overline{\B}) @>>>  H^0_J(C) @>>> H^1_J(\image(\varphi))
\cong H^1_J(\overline{\A}/xt\overline{\A}) \, . \end{CD}
\] Therefore $C=H^0_J(C) \hookrightarrow
H^1_J(\overline{\A}/xt\overline{\A})$.

By Lemma~\ref{S}(c), we have $H^1_J(\overline{\A})=0$. Hence 
$0 \rightarrow xt\overline{\A} \rightarrow
\overline{\A} \rightarrow \overline{\A}/xt \overline{\A}
 \rightarrow 0$ gives the exact sequence
\[\begin{CD} 0=H^1_J(\overline{\A}) @>>> H^1_J(\overline{\A}/xt\overline{\A})
@>>> H^2_J(xt\overline{\A}) \, .\end{CD}
\] Therefore $C \hookrightarrow H^2_J(xt \overline{\A})$.

Let $L={\rm ann \,}_{\overline{\A}}(xt)$. Notice that $xt$ is a generic
element of $I't$. Therefore $L=H^0_{I't}(L)$, which implies
$H^i_{I't}(L)=0$ for all $i \geq 1$. It follows that $H^i_{I't}(L_{t^{-1}})=0$ for
all $i \geq 1$. Since $J=(t^{-1},I't)$, we have an exact sequence
(\cite[8.1.2]{BS})
\[
H^1_{I't}(L_{t^{-1}}) \longrightarrow H^2_{J}(L)
\longrightarrow H^2_{I't}(L) \longrightarrow
H^2_{I't}(L_{t^{-1}})\longrightarrow H^3_{J}(L)
\longrightarrow H^3_{I't}(L) \, .
\] Therefore $H^2_{J}(L)=H^3_{J}(L)=0$. The exact sequence $0 \rightarrow L \rightarrow
\overline{\A}(-1) \stackrel{xt}{\longrightarrow} xt \overline{\A}
\rightarrow 0$ now yields an exact sequence
\[\begin{CD} 0=H^2_J(L) @>>> H^2_J(\overline{\A}(-1)) @>>> H^2_J(xt \overline{\A}) @>>> H^3_J(L)=0 \, , \end{CD} \]
from which we conclude that $C_n \hookrightarrow \left[H^2_J(xt
\overline{\A})\right]_n \cong \left[H^2_J(\overline{\A})
\right]_{n-1}$. By Theorem~\ref{H2neg}, we have
$\left[H^2_J(\overline{\A}) \right]_{n-1}=0$ for all $n-1 \leq 0$.
In particular $C_1=0$, that is, $\overline{I'/(x)} =
\overline{I'}/(x)$.
\end{proof}

\smallskip

\begin{Corollary}\label{ext1}
Let $R$ be a Noetherian, locally equidimensional, universally
catenary ring such that $R_{\red}$ is locally analytically
unramified. Let $I=(a_1, \ldots, a_n)$ be an $R$--ideal of height at
least $2$. Let $R'=R[z_1, \ldots, z_n]$, $I'=IR'$, and $x={
\sum_{i=1}^{n}} z_ia_i$. Then $\sqrt{(x)} \subset \overline{I'}$.
\end{Corollary}

\begin{proof} Let $y \in R'$ such that $y^m \in (x)$ for some $m$. Then
$(y+(x))^m =0$ in $R'/(x)$. In particular, $y+(x) \in
\overline{I'/(x)}$. Since $\overline{I'/(x)}=\overline{I'}/(x)$ by
Theorem~\ref{extension}, we get $y \in \overline{I'}$. 
\end{proof}

The reader may want to compare the previous corollary to the known result that
if $R$ is a Noetherian reduced ring, $I$ an $R$--ideal of grade at least $2$,
and $x$ a generic element for $I$ as above, then  $\sqrt{(x)}=(x)$ (\cite[Theorem (b)]{Hoc}).

The next corollary is a variation and an immediate consequence of Theorem \ref{extension}. We assume
$(R, \mathfrak{m})$ is local and we replace the polynomial ring $R'=R[z_1, \ldots, z_n]$
by its localization $R''=R(z_1, \ldots, z_n) = R[z_1, \ldots, z_n]_{\mathfrak{m}R[z_1, \ldots,
z_n]}$. In this case, the rings $R$ and $R''$ have the same dimension and the
generic element $x$ is part of a minimal generating set of the ideal $IR''$.

\begin{Corollary}\label{ext2}
Let $R$ be a Noetherian, equidimensional,
universally catenary local ring such that $R_{\red}$ is analytically
unramified. Let $I=(a_1, \ldots, a_n)$ be  an $R$--ideal of height at
least $2$. Let $R''=R(z_1, \ldots, z_n)$, $I''=IR''$, and $x={ \sum_{i=1}^{n}} z_ia_i$. Then $
\overline{I''/(x)} = \overline{I''}/(x)$.
\end{Corollary}

So far we have shown that the integral closure of an ideal is preserved 
by specialization modulo a generic element. One of the main
applications of this result is in proofs using induction on
the height of an ideal. For example, under slightly modified
assumptions, it yields a direct proof of a
well-known theorem on integral closures of ideal powers, due independently to 
Huneke and Itoh (\cite[4.7]{H1}, \cite[Theorem 1]{Itoh}).

\begin{Theorem}\label{Itoh}
Let $R$ be a Noetherian, locally equidimensional, universally
catenary ring such that $R_{\red}$ is locally analytically
unramified. Let $I$ be a complete intersection $R$--ideal. Then
\[ \overline{I^{n+1}} \bigcap I^n = \overline{I}I^n \quad \mbox{\rm for all} \;\; n \geq 0 \, .\]
\end{Theorem}

\begin{proof} We may assume that $R$ is local and $g=\height(I) >0$.
Write $ G = { \bigoplus_{n \geq 0}} I^n/I^{n+1} $ for the associated graded ring of $I$ and let
$\widetilde{G}={\bigoplus_{n \geq 0}} \overline{I^n}/\overline{I^{n+1}}$.
Consider the exact sequence
\[\begin{CD} 0 @>>> K @>>> G \otimes_R R/\overline{I} @>{\varphi}>>
\widetilde{G} \, , \end{CD} \] where $\varphi$ is the natural map and
$K$ is its kernel. Notice that
the degree $n$ component of $K$ is $K_n={\ds \frac{\overline{I^{n+1}}
\bigcap I^n}{\overline{I}I^n}}$. We want to show that $K=0$.

We use induction on $g$. First let $g=1$. In this case $I$ is
generated by a single regular element $a$. Write $\overline{R}$ for the
integral closure of $R$ in its total ring of quotients. Since
$\overline{I}=a \overline{R} \cap R$ and
$\overline{I^{n+1}}=a^{n+1}\overline{R} \cap R$, we obtain
$$\overline{I^{n+1}} \cap I^n= a^{n+1}\overline{R} \cap a^nR = a^n \left( a \overline{R} \cap R \right) =\overline{I}I^n \, .$$ 
It follows that $K=0$ in this case.

Suppose that $g \geq 2$. Let $I=(a_1, \ldots, a_g)$,
$R''= R(z_1, \ldots, z_g)$, and $I''=IR''$. 
Denote the associated graded ring of $I''$ by $G''$, set $\widetilde{G''}={ \bigoplus_{n \geq 0}}
\overline{(I'')^n}/\overline{(I'')^{n+1}}$, and consider the natural exact sequence
\[\begin{CD} 0 @>>> K'' @>>> G'' \otimes_{R''} R''/\overline{I''} @>{{\varphi}''}>>
\widetilde{G''} \, . \end{CD} \]
Further we define $x=\sum_{i=1}^{g} z_ia_i
\in I''$, $\mathfrak{R}=R''/(x)$, and $\mathfrak{I}=I''/(x)$.
We wish to apply the induction hypothesis to the ideal $\mathfrak{I}$
of the ring $\mathfrak{R}$. Notice that $\mathfrak{I}$ is a complete
intersection ideal of height $g-1$.

We first argue that the condition concerning analytic unramifiedness
passes from $R$ to $\mathfrak{R}$. For this part we may assume that
$\widehat{R}$ is reduced and it suffices to show that $\widehat{\mathfrak{R}}$ is reduced as well.
But indeed, $\widehat{\mathfrak{R}}$ is the completion of the excellent local ring $\widehat{R} (z_1, \ldots, z_g)/ (x)$,
and the latter ring is reduced because $x$ is a generic element
for $\widehat{I}$, an ideal of grade at least $2$ in the reduced ring $\widehat{R}$ (\cite[Theorem (b)]{Hoc}).

By Corollary~\ref{ext2}, we have
$\overline{\mathfrak{I}}=\overline{I''}/(x)$. Let $\mathcal{G}$ be the associated
graded ring of $\mathfrak{I}$. We observe that $\mathcal{G} \cong G''/(x+(I'')^2)$ because $x$
is part of a regular sequence that generates the ideal $I''$. Moreover, the element
$x+(I'')^2$ is regular on $\widetilde{G''}$, hence on ${\rm Im}({\varphi}'')$
(Lemma~\ref{S}(c)). Write $\widetilde{\mathcal{G}}={\bigoplus_{n \geq 0}}
\overline{\mathfrak{I}^n}/\overline{\mathfrak{I}^{n+1}}$. We
obtain the following commutative diagram with exact rows,
\[
\begin{CD}
0 @>>> \mathcal{G} \otimes_{G''}K''  @>>> \mathcal{G} \otimes_{G''} G'' \otimes_{R''} R''/\overline{I''}
@>>> \mathcal{G} \otimes_{G''}{\rm Im}({\varphi}'') @>>> 0 \\
&& && @V{\cong}VV @VVV \\
&&0 @>>> \mathcal{G} \otimes_{\mathfrak{R}}
\mathfrak{R}/\overline{\mathfrak{I}} @>>> \widetilde{\mathcal{G}} \ ,
\end{CD}
\]
where the vertical isomorphism follows from the equality
$\overline{\mathfrak{I}}=\overline{I''}/(x)$, the exactness of the
first row from the fact that $x+(I'')^2$ is regular on  ${\rm Im}({\varphi}'')$, and the exactness of the second row from the
induction hypothesis. We conclude that $\mathcal{G}\otimes_{G''}K''=0$, which gives $K''=0$ 
by the graded Nakayama Lemma. Therefore $K=0$, as asserted.
\end{proof}

\bigskip

\section{Modules over two dimensional regular local rings}

\smallskip

\begin{Notation}\label{notation}{\rm
Let $(R, \mathfrak{m})$ be a Noetherian local ring and $E$ a
finitely generated torsion-free $R$--module such that $E_{\mathfrak
p}$ is free for every prime ideal $\mathfrak p$ of $R$ with $\depth
R_{\mathfrak p} \leq 1$. In this case $E$ has a rank
$e$, which we assume to be positive (\cite[2.1]{Hartshorne}). 
Let $U=\sum_{j=1}^n Ra_j$ be a submodule of $E$ so that $E/U$ has 
grade at least $2$, in other words, $U_{\mathfrak p}=E_{\mathfrak p}$ 
whenever $\depth R_{\mathfrak p} \leq 1$. This condition is automatically
satisfied if $U$ is a reduction of $E$ because in that case $E/U$ is 
supported on the non-free locus of $E$.
Let $Z=\{z_{ij}\;|\; 1 \leq i \leq e-1,\; 1 \leq j
\leq n\}$ be a set of indeterminates over $R$. Set
\[
R''=R[Z]_{\mathfrak{m}R[Z]},\;\; U''=R'' \otimes_R U,\;\; E ''=R''
\otimes_R E,\;\; x_i=\sum_{j=1}^{n} z_{ij}a_{j}, \;\; \mbox{\rm and}\;
F=\sum_{i=1}^{e-1} R''x_{i} \, .
\] 

The module $F$ is a free $R''$--module of rank $e-1$ and $E''/F$ is
torsion-free with rank $1$, hence isomorphic to an $R''$--ideal
(\cite[3.2]{SUV}). An $R''$--ideal $I$ with $I \cong E''/F$ is
called a {\em generic Bourbaki ideal} of $E$ with respect to $U$
(\cite[3.3]{SUV}).  We refer the reader to \cite{SUV} for a general
discussion of generic Bourbaki ideals, including the fact that a
generic Bourbaki ideal is essentially unique (\cite[3.4]{SUV}). If
$E$ has finite projective dimension, then the ideal $I \cong E''/F$
can be chosen to have grade at least $2$ (\cite[3.2]{SUV}). This
ideal is uniquely determined by $E$ and the generators $a_1, \ldots,
a_n$ of $U$; we will call it a {\em generic Bourbaki ideal} of $E$
with respect to $U$ {\em of grade at least $2$}. If $U=E$, we simply
call it a generic Bourbaki ideal of $E$ of grade at least $2$.}
\end{Notation}

\medskip

Our objective is to study whether a generic Bourbaki ideal of the
integral closure $\overline{E}$ with respect to $E$ is integrally
closed. In this section we consider the case of two-dimensional
regular local rings. We prove that a generic Bourbaki ideal of
grade at least $2$ of an integrally closed module is again 
integrally closed and we express this ideal as a Fitting ideal 
of the module. Using these facts, we deduce the known results about 
integrally closed modules over two-dimensional regular local 
rings from the corresponding classical theory for ideals. 

\smallskip

Thus let $(R, \mathfrak{m})$ be a two-dimensional regular local ring.
For an $R$--ideal $I$, the {\em order} $o(I)$ of $I$ is defined to
be the largest integer $r$ such that $I \subset \mathfrak{m}^r$.
Let $E$ be a finitely generated torsion-free $R$--module with rank
$e$ and write $E^*=\Hom_R(E,R)$. There is a natural embedding $E
\subset E^{**}$, where $E^{**}$ is a free module of rank $e$ and $E^{**}/E$ 
has finite length. One defines the {\em
order} $o(E)$ of $E$ to be the order of the zeroth Fitting 
ideal $\mbox{\rm Fitt}_0(E^{**}/E)$. We recall that for a presentation 
$R^m \stackrel{\varphi}{\longrightarrow} R^n \rightarrow M \rightarrow 0$ 
of a module $M$, the $i$-{\it th Fitting ideal} $ \, \mbox{\rm Fitt}_i(M)$ of 
$M$ is the $R$--ideal $I_{n-i}(\varphi)$ generated by the $(n-i) \times
(n-i)$-minors of $\varphi$. This ideal only depends on $M$ and $i$.

Suppose that the maximal ideal $\mathfrak{m}$ is
generated by $a,b$ and let $S=R\left[\frac{\mathfrak{m}}{a}\right]$.
Identifying $E^{**}$ with $R^e$, which is contained in $S^e$, we let
$ES$ denote the $S$--submodule of $S^e$ generated by $E$. 
This module 
$ES$ is called the {\em transform} of $E$ in $S$. It 
is isomorphic to $E \otimes_R S$ modulo $S$--torsion. The
$R$--module $E$ is said to be {\em contracted} if $E=ES \cap R^e$ for
some choice of $a$. Finally, we denote the minimal number of generators of $E$
by $\nu(E)$.

We will use these facts about transforms of modules and contracted
modules that are proved in \cite{K}:

\begin{Theorem}\label{transform}
Let $R$ be a $2$-dimensional regular local ring with  maximal ideal
$\mathfrak{m}=(a,b)$ and let $E$ be a finitely generated torsion-free
$R$--module.
\begin{enumerate}
\item[{\rm (a)}] If $E$ is integrally closed, then so is the transform
$ES$ in $S=R\left[\frac{\mathfrak{m}}{a}\right]$.

\item[{\rm (b)}] If $E$ is integrally closed, then $E$ is
contracted.
\item[{\rm (c)}] The module $E$ is contracted if and only if
$\nu(E)=o(E)+ \rank(E)$.
\end{enumerate} 
\end{Theorem}

\begin{proof}  Part (a) is proved in \cite[4.6]{K},  part (b) in
\cite[4.3]{K}, and  part (c) in \cite[2.5]{K}. 
\end{proof}

\vspace{.0000000001in}

\begin{Lemma}\label{Fitt}
Let $R$ be a $2$-dimensional regular local ring.
Let $E$ be a finitely generated torsion-free $R$--module with rank $e
>0$ and let $U$ be a submodule of $E$ so that $E/U$ has finite length. 
Let $I \cong
E''/F$ be a generic Bourbaki ideal of $E$ with respect to $U$ of
grade at least $2$. Then $I= \,\mbox{\rm Fitt}_1(E''/F)$ and
this ideal is a reduction of $\,\mbox{\rm Fitt}_e(E) R''$.
\end{Lemma}

\begin{proof} We may assume that the residue field of $R$ is infinite.
Consider an epimorphism $R''^n \twoheadrightarrow E''$ and
decompose $R''^n \cong F \oplus G$, where $G$ is a free $R''$--module generated
by $n-e+1$ general elements of $R''^n$. Let 
$0 \rightarrow R''^{n-e} \stackrel{\varphi}{\longrightarrow} F \oplus G \longrightarrow
E'' \rightarrow 0$ be a presentation of $E''$, and let $\psi$ be the $n-e+1$ by $n-e$ submatrix of
$\varphi$ consisting of the last rows of $\varphi$. Then $\psi$ is a
presentation matrix for the module $E''/F \cong I$. Since $I$ is
an ideal of projective dimension at most $1$ and grade at least $2$, the Hilbert-Burch Theorem
gives that 
$I=I_{n-e}(\psi)$. Clearly $\Fitt_1(E''/F)=I_{n-e}(\psi) \subset I_{n-e}(\varphi)= \Fitt_e(E)R''$.

It remains to prove that $I_{n-e}(\varphi)$ is integral over $I_{n-e}(\psi)$. Let $-^*$
denote dualizing into $R''$. We consider the map
$F^* \oplus G^* \stackrel{{\varphi}^*}{\longrightarrow} (R''^{n-e})^*$. Its image
${\rm{Im}}({\varphi}^*)$ has rank $n-e$ as a module over $R''$, a $2$-dimensional 
local ring with infinite residue field. Therefore $(n-e)+2-1 =n-e+1$ general elements
of this module generate a reduction (\cite[2.3]{SUV}). Thus ${\varphi}^*(G^*)$ 
is a reduction of ${\rm{Im}}({\varphi}^*)$. In other words, the colums of the matrix
${\varphi}^*$ are integral over the module generated by the colums of ${\psi}^*$.
Now the valuative criterion for integral dependence shows that 
$I_{n-e}({\varphi}^*)$ is integral over $I_{n-e}({\psi}^*)$. Thus indeed
$I_{n-e}(\varphi)$ is integral over $I_{n-e}(\psi).$
\end{proof}

\vspace{.0000000001in}

\begin{Lemma}\label{o(E)}
Let $R$ be a $2$-dimensional regular local ring.
Let $E$ be a finitely generated torsion-free $R$--module with rank $e
>0$ and let $U$ be a submodule of $E$ so that $E/U$ has finite length. Let $I \cong
E''/F$ be a generic Bourbaki ideal of $E$ with respect to $U$ of
grade at least $2$.

\begin{enumerate}
\item[{\rm (a)}] The Fitting ideals $\,\mbox{\rm Fitt}_0(E^{**}/E)$
and $\,\mbox{\rm Fitt}_e(E)$ are equal. Hence
$o(E)=o(\mbox{\rm Fitt}_e(E))$.
\item[{\rm (b)}] The order of $E$ is equal to the order of
$I$.
\item[{\rm (c)}] If $E$ is contracted, then I is contracted.
\end{enumerate}
\end{Lemma}

\begin{proof}  (a) Since $E^{**}$ is free and $E^{**}/E$ has projective dimension at most $2$,  the
ideals $\mbox{\rm Fitt}_0(E^{**}/E)$ and $\mbox{\rm Fitt}_e(E)$ are
isomorphic (\cite[3.1]{BE}). But both have grade at least $2$,
which implies that they are equal.

\smallskip

\noindent (b) One has $o(E)=o({\rm Fitt}_e(E))=o({\rm Fitt}_e(E)R'')=o(I)$. 
The first equality holds by part (a) and the second equality is obvious.
The third equality follows from Lemma \ref{Fitt}; indeed, the lemma shows that
${\rm Fitt}_e(E)R''$ is integral over $I$, and any ideal of a regular 
local ring has the same order as its integral closure because the powers of the 
maximal ideal are integrally closed.

\smallskip

\noindent (c) Suppose $E$ is contracted. Then $E''$ is contracted so
that $\nu(E'')=o(E'')+e$ by Theorem~\ref{transform}(c). Hence we
obtain
\[
\nu( I) \leq o(I) + 1 =
o(E'')+1=\nu(E'')-e+1 \leq \nu(I),
\]where the first inequality holds by the Hilbert-Burch Theorem, the first equality follows from part (b), the
second equality is explained above, and the last inequality is obvious from the isomorphism $I \cong E''/F$.
Therefore $I$ is contracted according to Theorem~\ref{transform}(c).  
\end{proof}

\medskip

Let $R$ be a $2$-dimensional regular local ring with maximal
ideal $\mathfrak{m}=(a,b)$ and write
$S=R\left[\frac{\mathfrak{m}}{a}\right]$. Let $\mathfrak{N}$ be a
maximal ideal of $S$ containing $\mathfrak{m}S$. Then
$T=S_{\mathfrak{N}}$ is a $2$-dimensional regular local ring and is
called a {\em first quadratic transformation} of $R$.  For an
$R$--ideal $I$ of grade at least $2$, we can write $IS=a^r I^{\sharp} $ for some $S$--ideal
$I^{\sharp}$, where $r=o(I)$. We call the localization
$(I^{\sharp})_{\mathfrak{N}}$ a {\em first quadratic transform} of
$I$ and denote it by $\widetilde{I}$. This ideal has grade at least $2$
(\cite[14.3.2]{HSw}). Let $E$ be a finitely generated torsion-free
$R$--module. We call the module $ET= (ES)_{\mathfrak{N}}$ the {\em transform} of $E$ in
$T$. The next lemma explains how a generic Bourbaki ideal of grade at least $2$
behaves under a first quadratic transformation.

\begin{Lemma}\label{BItransform} Let $(R, \m)$ be a $2$-dimensional
regular local ring and let $\, T = R \left[ \frac{\m}{a} \right]_{\mathfrak{N} }$ 
be a first quadratic transformation of $R$. Let $E$ be a finitely generated torsion-free
$R$--module and $U$ a reduction of $E$.  
Let $I \subset R''$ be a generic Bourbaki ideal of $E$ with respect to 
$U$ of grade at least $2$ and write 
$T'' = R'' \left[ \frac{\m}{a} \right]_{\mathfrak{N} R'' \left[ \frac{\m}{a} \right] }$. 
Then the first quadratic transform of $I$ in $T''$ 
is a generic Bourbaki ideal of $ET$ with respect to $UT$
of grade at least $2$.
\end{Lemma}

\begin{proof} We use the notation of \ref{notation}; in particular, $R''=R[Z]_{{\mathfrak m}R[Z]}$ and $I \cong E''/F$. 
Notice that $T''=T[Z]_{\mathfrak{N}T[Z]}$, the module $UT$ is a reduction of $ET$, and
$FT''$ is a generic $(e-1)$-generated submodule of $UT''$. Therefore $E''T''/FT''$ is a torsion-free
$T''$--module of rank $1$. There are surjective homomorphisms of $T''$--modules of rank $1$,
\[ E''T''/FT'' {\twoheadleftarrow} \left( E''/F \right) \otimes_{R''} T'' {\twoheadrightarrow} IT''. \]
Since the outer two modules are torsion-free, one obtains an isomorphism
$E''T''/FT'' \cong IT''$. As $IT'' \cong 
\widetilde{I}$ and the latter ideal has grade at least $2$,
it follows that $\widetilde{I}$ is a generic Bourbaki ideal of $ET$ with respect to $UT$ of grade at
least $2$.
\end{proof}

\smallskip

\begin{Theorem}\label{2dimRLR}
Let $R$ be a $2$-dimensional regular local ring.
Let $E$ be a finitely generated torsion-free $R$--module with rank $e > 0$ and let $U$ be a reduction of $E$.
Then $E$ is integrally closed if and only if any generic 
Bourbaki ideal of $E$ with respect to $U$ is integrally closed.
\end{Theorem}

\begin{proof} It will follow from Proposition~\ref{HImodconv}
that $E$ is integrally closed if some generic 
Bourbaki ideal of $E$ with respect to $U$ is integrally closed.
To prove the converse suppose that $E$ is integrally closed.
We may assume that the residue field of $R$ is infinite.
We use the notation of \ref{notation} and we write
$\mathfrak{n}= \mathfrak{m} R''=(a,b)$ for the maximal ideal of $R''$. 
Furthermore, let $I \cong E''/F$ be a generic Bourbaki ideal of $E$ with 
respect to $U$ of grade at least $2$. 

In order to prove that any
generic Bourbaki ideal of $E$ with respect to $U$ is integrally
closed, it suffices to show that $I$ is integrally closed (\cite[3.4(a)]{SUV}). 
We use induction on the multiplicity
$e(I)$ of $I$. If $e(I)=0$, then $I=R''$ and we are done. Suppose
$e(I) \geq 1$. Since $E$ is integrally closed, $E$ is contracted
according to Theorem~\ref{transform}(b). Therefore $I$ is contracted
by Lemma~\ref{o(E)}(c). Hence
we have $I=IS'' \cap R''$ for some $S''=R''\left[
\frac{\mathfrak{n}}{a}\right]$. Since $R/\m$ is infinite, we may choose 
$a$ to be a general element of $\m$ (\cite[14.2.2]{HSw}), in particular, 
$a \in \m$ and $S''=R''[\frac{\mathfrak{m}}{a}]$. As $I$ is contracted from 
$IS''$, it suffices to 
prove that $IS''$ is integrally closed. Write 
$I^{\sharp}=\frac{1}{a^r}IS''$, where $r=o(I)$. We need to show that 
$(I^{\sharp})_{\mathfrak{N''}}$ is integrally closed for every maximal 
ideal $\mathfrak{N''}$ containing $I^{\sharp}$. Notice that any such $\mathfrak{N''}$ 
contains $\mathfrak{m}$. For simplicity, write
$T''=S''_{\mathfrak{N''}}$, $\widetilde{I}=(I^{\sharp})_{\mathfrak{N''}}$, 
$S=R \left[ \frac{\m}{a} \right]$, ${\mathfrak{N}}={\mathfrak N''} \cap S$,
and $T=S_{\mathfrak{N}}$. 

First assume that ${\dim} \,  T  \geq 2 = {\dim} \, T'' $. Since $S''$ is a 
ring of fractions of the polynomial ring $S[Z]$, it follows that 
$\mathfrak{N''}$ is extended from $\mathfrak{N}$ in this case, hence 
$T''=S''_{{\mathfrak N}S''}$. Thus Lemma~\ref{BItransform} shows 
that $\widetilde{I}$ is a generic Bourbaki ideal of $ET$ with 
respect to $UT$ of grade at least $2$. The transform $ET$ of the 
integrally closed module $E$ is again integrally 
closed according to Theorem~\ref{transform}(a).
Moreover, \cite[3.6]{H} gives $e(\widetilde{I}) < e(I)$.
Now our induction hypothesis implies that $\widetilde{I}$
is integrally closed, as asserted.

Next suppose that ${\dim} \, T  \leq 1 $. We will show that the $T''$--module
$ET''/FT''$ is cyclic in this case. Since this module maps onto $\widetilde{I}$,
it will follow that $\widetilde{I}$ is principal, hence integrally closed.
Now, since ${\mathfrak m} \subset {\mathfrak N}$, we see that ${\dim} \,  T  = 1$. 
Thus $T$ is a discrete valuation ring. In particular, the local map $T \to T''=T[Z]_{{\mathfrak M}}$
is flat with a one-dimensional closed fiber $K[Z]_{{\mathfrak M}K[Z]}$, 
where ${\mathfrak M}$ is a prime ideal of $T[Z]$ and $K$ stands for
the residue field of $T$. Moreover, the $T$--module $ET$ is free
of rank $e$ and therefore $UT=ET$. After changing the generating sequence $a_1, \ldots, a_n$ of the module
$UT=ET$ and applying a linear change of variables of the polynomial ring $T[Z]$, we
may assume that $a_1, \ldots, a_e$ form a basis of $ET$ and $a_{e+1}= \ldots =a_n=0$.
Thus ${\rm Fitt}_1(ET''/FT'')=I_{e-1}(\varphi)T''$, where $\varphi$ is the $e-1$ by $e$ matrix
consisting of the first $e$ columns of the matrix of variables $\left( z_{ij} \right)$.
This ideal has height at least two when extended to $K[Z]_{{\mathfrak M}K[Z]}$, a ring of
dimension one. We conclude that the extended ideal is the unit ideal and hence so is ${\rm Fitt}_1(ET''/FT'')$.
Thus the $T''$--module $ET''/FT''$ is cyclic, as asserted.
\end{proof}

\vspace{.0000000001in}

\begin{Corollary}\label{Fe(E)}
Let $R$ be a $2$-dimensional regular local ring.
Let $E$ be a finitely generated torsion-free $R$--module with rank $e > 0$ and 
let $U$ be a reduction of $E$. Let $I \subset R''$ be a
generic Bourbaki ideal of $E$ with respect to $U$ of grade at least $2$. 
If $E$ is integrally closed, then $I=\mbox{\rm Fitt}_e(E)R''$ and $ \, \mbox{\rm Fitt}_e(E)$
is integrally closed.
\end{Corollary}

\begin{proof} The assertion follows from Lemma~\ref{Fitt} and
Theorem~\ref{2dimRLR}.
\end{proof}

That the ideal $\mbox{\rm Fitt}_e(E)$ in the corollary above is 
integrally closed can also be seen from \cite[5.4]{K} by
way of Lemma~\ref{o(E)}(a).
The next result about integrally closed modules over
$2$-dimensional regular local rings has been proved in
\cite[5.2]{K} and \cite[4.1]{KK}. Here we deduce it from the
corresponding statements for ideals, using Theorem~\ref{2dimRLR}.

\begin{Corollary}\label{application1} Let $R$ be a $2$-dimensional
regular local ring. Let $E$ be a finitely generated torsion-free
$R$--module. Suppose $E$ is integrally closed. Then
\begin{enumerate}
\item[{\rm (a)}] The module $E$ is normal.
\item[{\rm (b)}] $E^2=UE$ for every reduction $U$ of $E$.
\item[{\rm (c)}] The Rees algebra $\R(E)$ is Cohen-Macaulay.
\item[{\rm (d)}] Let $G$ be a finitely generated torsion-free integrally closed $R$--module. Then $EG=(E\otimes_R G)/\tau(E
\otimes_R G)$ is integrally closed, where $\tau(-)$ denotes
$R$--torsion.
\end{enumerate}
\end{Corollary}

\begin{proof} Let $U$ be any reduction of $E$. By Theorem~\ref{2dimRLR}, 
a generic Bourbaki ideal of $E$ with respect to $U$ of grade
at least $2$ is integrally closed. Now  part (a) follows from
\cite[14.4.4]{HSw} and \cite[3.5(a)(ii)]{SUV}, part (b) from
\cite[5.5]{LT} and \cite[3.5(b),(c)]{SUV}, and part (c) from
\cite[3.2]{HS} and \cite[3.5(a)(i)]{SUV}. Finally for part (d),
notice that $(E \oplus G)^2$ is integrally closed by part (a). Since
$EG$ is a direct summand of $(E \oplus G)^2=E^2 \oplus EG \oplus
G^2$, the module $EG$ is integrally closed as well. 
\end{proof}

\medskip

\section{Specialization by generic elements of modules}

\smallskip

In this section we are going to show, under fairly general assumptions,
that integral closures of modules are compatible with factoring out generic
elements. Some of the proofs here will rely on the results for 2-dimensional 
regular local rings that were obtained in the previous section. 

\begin{Lemma}\label{emb0}
Let $R$ be a Noetherian ring and $E$ a finitely generated
torsion-free $R$--module having a rank $e$. Suppose that $E_{\mathfrak p}$ is
free whenever $\depth R_{\mathfrak{p}} \leq 1$. Then there exists an embedding $E \subset
R^e$ such that $(R^e/E)_{\mathfrak{p}}$ is cyclic whenever $\depth
R_{\mathfrak{p}} \leq 1$.
\end{Lemma}

\begin{proof} We write $-^* = \Hom_R(-,R)$. There is an epimorphism
$(R^m)^* \twoheadrightarrow E^*$, which induces an embedding $E^{**}
\hookrightarrow R^m$. Since the natural homomorphism $E
\rightarrow E^{**}$ is injective, we obtain an exact sequence $0
\rightarrow E \rightarrow R^m \rightarrow C \rightarrow 0$.
For every prime ideal $\mathfrak p$ with $\depth R_{\mathfrak{p}} \leq 1$, 
the $R_{\mathfrak p}$--modules $E_{\mathfrak{p}}$ and $E^*_{\mathfrak{p}}$ are free of rank $e$, 
hence $C_{\mathfrak{p}}$ is free of rank $m-e$. By basic element theory there 
exists a basis $\alpha_1, \ldots , \alpha_m$ of $R^m$ such that for the $R$--modules
$G \subset F$ generated by the images in $C$ of  $\alpha_1, \ldots , \alpha_{m-e-1}$ and
of  $\alpha_1, \ldots , \alpha_{m-e}$, 
the $R$--module $C/F$ has rank $0$ and $(C/G)_{\mathfrak{p}} \cong R_{\mathfrak{p}}$ whenever
$\depth R_{\mathfrak{p}} \leq 1$
(\cite[1.1]{MY}). In particular, $F$ is a free $R$--module of rank $m-e$. Thus
we obtain the following commutative diagram with exact rows and columns,
\[\begin{CD}
&&  &&   0 && 0\\
&&  &&   @VVV @VVV\\
&& 0 @>>> F @>>> F @>>> 0 \\
&&  @VVV @VVV @VVV\\
0 @>>> E @>>> R^m @>>> C @>>> 0 \ .
\end{CD} 
\]

\bigskip

\noindent
Since $F$ is a free direct summand of $R^m$ of rank $m-e$, the 
Snake Lemma yields an exact sequence 
$ 0 \rightarrow E \rightarrow R^{e} \rightarrow C/F \rightarrow
0$.
For any prime ideal $\mathfrak{p}$ with $\depth R_{\mathfrak{p}} \leq 1 $, we have
\[\nu((R^e/E)_{\mathfrak{p}})=\nu((C/F)_{\mathfrak{p}}) \leq
\nu((C/G)_{\mathfrak{p}})    =1 \, .\] \end{proof}

Let $R$ be a Noetherian ring and $E$ a finitely generated $R$--module having a
rank. Recall that $E$ is said to be of {\it linear type} if the
natural map from the symmetric algebra $\S(E)$ onto the Rees algebra
$\R(E)$ is an isomorphism, or equivalently, if $\S(E)$ is
$R$--torsion-free.

\begin{Lemma}\label{emb1}
Let $R$ be a Noetherian ring and let $E={ \sum_{i=1}^n} Ra_i$ be an
$R$--module having a rank $e \geq 2$. Let $E \subset
R^e$ be an embedding such that $\nu((R^e/E)_{\mathfrak{p}}) \leq
e-1$ whenever $\depth
R_{\mathfrak{p}} \leq 1$. Let $R'=R[z_1, \ldots, z_n]$, $E'=R' \otimes_{R} E$, and
$x={ \sum_{i=1}^n} z_ia_i$. Then $R'^e/R'x$ is of linear type.
\end{Lemma}

\begin{proof} The embedding $E \subset R^e$ and the chosen generators 
$a_1, \ldots, a_n$ of $E$ induce an
exact sequence $ R^n \stackrel{\psi}{\longrightarrow} R^e
\rightarrow R^e/E \rightarrow 0$. Notice that $\grade I_1(\psi) \geq 2$
because $\nu((R^e/E)_{\mathfrak{p}}) \leq e-1$ for every $\mathfrak{p} \in
\Spec(R)$ with $\depth R_{\mathfrak{p}} \leq 1$. Consider the column vector
$\Phi=\psi \left[\begin{array}{c} z_1 \\ \vdots \\
z_n
\end{array}  \right]$.
Since $\grade I_1(\Phi) \geq 1$ 
we have an exact sequence $0 \rightarrow R'
\stackrel{\Phi}{\longrightarrow} R'^e \longrightarrow R'^e/R'x
\rightarrow 0$. As $R'^e/R'x$ has projective dimension at
most $1$, in order to prove that this module is of linear type, it suffices
to show that $\grade I_1(\Phi) \geq 2$ (\cite[Proposition 4]{A}).

Suppose there exists a prime ideal $\mathfrak{P} \in \Spec(R')$
containing $I_1(\Phi)$ such that $\depth R'_{\mathfrak{P}} \leq 1$.
Let $\mathfrak{p} = \mathfrak{P} \cap R$. We have $\depth
R_{\mathfrak{p}} \leq 1$. Then since $\grade I_1(\psi) \geq 2$, we get
$I_1(\psi)_{\mathfrak{p}}={R}_{\mathfrak{p}}$. After elementary row and 
column operations on the matrix $\psi_{\mathfrak p}$ and a linear change
of variables of the polynomial ring $R_{\mathfrak p}[z_1, \ldots, z_n]$, we may assume that
\[
\psi_{\mathfrak{p}} = \left[\begin{array}{cccc} 1 &  0 & \cdots & 0 \\
0  & b_{22} &
\cdots & b_{2n} \\
0  & \vdots &
\cdots & \vdots \\
0  & b_{e2} & \cdots & b_{en}
\end{array} \right] \quad \mbox{\rm and} \quad
\Phi_{\mathfrak{P}} = \left[\begin{array}{c} z_1 \\ b_{22}z_2+ \cdots +b_{2n}z_n \\
\vdots \\ b_{e2}z_2+ \cdots + b_{en}z_n
\end{array} \right].
\]Since
$E={\rm Im}(\psi)$ has rank $e \geq 2$, we have $\grade I_2(\psi) \geq
1$. Therefore the ideal 
\[(b_{22}z_2+  \cdots  +b_{2n}z_n,\  \ldots  \ , b_{e2}z_2+ \cdots + b_{en}z_n)\] 
in $R_{\mathfrak{p}}[z_2, \ldots, z_n]$ has grade
at least $1$. It follows that the $R'_{\mathfrak{P}}$--ideal
$I_1(\Phi_{\mathfrak{P}})$ has grade at least $2$, which is impossible since this 
ideal is contained in $\mathfrak{P}R'_{\mathfrak{P}}$ and $\depth R'_{\mathfrak{P}} \leq 1$.
\end{proof}

\smallskip
Lemma~\ref{emb1} and the next Theorem~\ref{HImod1} require an embedding $E \subset R^e$ such that
$\nu((R^e/E)_{\mathfrak{p}}) \leq e-1$ for any $\mathfrak{p} \in
\Spec(R)$ with $\depth R_{\mathfrak{p}} \leq 1$. Such an
embedding exists by
Lemma~\ref{emb0} if $e \geq 2$ and the assumptions of the lemma are satisfied. Now we 
are ready to prove the module analogue of Theorem~\ref{extension}.

\begin{Theorem}\label{HImod1}
Let $R$ be a Noetherian, locally equidimensional, universally catenary
ring such that $R_{\red}$ is locally analytically
unramified. Let $E={ \sum_{i=1}^n }R a_i$ be an
$R$--module having a rank $e \geq 2$. Let $E \subset R^e$ be an
embedding such that $\nu((R^e/E)_{\mathfrak{p}}) \leq e-1$ whenever
$\depth R_{\mathfrak{p}} \leq 1$. Let $R'=R[z_1, \ldots, z_n]$,
$E'=R' \otimes_{R} E$, and $x= { \sum_{i=1}^n }z_ia_i$. Let
$\overline{E'}$ and $\overline{E'/R'x}^{\;R'^e/R'x}$ denote the
integral closure of $E'$ in $R'^e$  and that of $E'/R'x$ in
$R'^e/R'x$. Then
\[
\overline{E'/R'x}^{\;R'^e/R'x} = \overline{E'}/R'x \, .
\]
\end{Theorem}

\begin{proof} We may assume that $(R, {\mathfrak m})$ is a local ring of dimension $d$. The
embedding $E \subset R^e$
identifies $E$ with a set of linear forms in 
the polynomial ring $S=\S(R^e) =R[t_1, \ldots, t_e]$. Let $I$ be the
$S$--ideal generated by $E$. The symmetric algebra $\S(R^e/E)$ of
$R^e/E$ is $S/I$, and from \cite[2.6]{HR} we obtain
\[\height(I)=  d+e - \dim(\S(R^e/E))= \min_{\mathfrak{p} \in \mbox{\rm \tiny Spec}(R)}
\left\{ \; e + \height(\mathfrak{p}) - \nu(( R^e/E)_{\mathfrak{p}})\;
\right\}.\] For any $\mathfrak{p} \in \Spec(R)$, we have

\[e +
\height(\mathfrak{p}) - \nu(( R^e/E)_{\mathfrak{p}}) \geq
\left\{\begin{array}{lll}  \height(\mathfrak{p}) &\geq 2 & \mbox{\rm
if}\; \height(\mathfrak{p}) \geq 2 \\&&\\ \height(\mathfrak{p})+1 &=
2 & \mbox{\rm if}\; \height(\mathfrak{p}) = 1 \\ &&\\
e & \geq 2 & \mbox{\rm if}\; \height(\mathfrak{p}) = 0 \, .
\end{array}
\right.
\] Therefore the $S$--ideal $I$ generated by $E$ has height at least $2$.
Let $S'=S[z_1, \ldots, z_n], I'=IS'$, and $ x=\sum_{i=1}^n{z_ia_i}$.
We have
$\overline{I'}/(x)=\overline{I'/(x)} \subset S'/(x)$
by Theorem~\ref{extension}. The theorem applies because  
$S_{\rm{red}}$ is locally analytically unramified; in fact, since $I$ 
is a homogeneous ideal, it suffices that the localization
$\left(S_{\rm{red}}\right)_{({\mathfrak m}, t_1, \ldots, t_e)}$ is
analytically unramified; but the completion of this ring is $\widehat{R_{{\rm red}}}[[t_1, \ldots, t_e]]$,
which is reduced.

The embedding $ E' \subset R'^e$ defines the Rees algebra $\R(E')$ as
a subalgebra of the polynomial ring $S'=\S(R'^e) =R'[t_1, \ldots,
t_e]$. Notice that $I'$ is the $S'$ ideal generated by $E'$. By
\cite[1.2]{HUV}, we have $\overline{E'}=\left[\; \overline{I'}
\;\right]_1$.
Therefore we obtain
\[
\left[ \; \overline{I'/(x) } \; \right]_1= \left[\;
\overline{I'}/(x) \; \right]_1=\overline{E'}/R'x \subset
\overline{E'/R'x}^{R'^e/R'x} \subset R'^e/R'x \, .\]

It remains to show that $\overline{E'/R'x}^{R'^e/R'x} \subset \left[
\; \overline{I'/(x) } \; \right]_1$. Since the module $R'^e/R'x$ is of
linear type (Lemma~\ref{emb1}), we have
\[
\R(R'^e/R'x) \cong \S(R'^e/R'x) \;\stackrel{\Phi}{\longrightarrow}\;
S'/(x) \, ,
\] where $\Phi$ is the natural map induced by the identity on $R'^e/R'x$.
Let $u \in \overline{E'/R'x}^{R'^e/R'x} \subset R'^e/R'x$. Then there is
an equation in $\R(R'^e/R'x)$ of the form
\begin{eqnarray*} u^m + a_1u^{m-1} + \cdots + a_m=0, \quad a_i\in
(E'/R'x)^i.\end{eqnarray*} Applying $\Phi$ to this equation, it
converts into an equation of integrality of $u \in S'/(x)$ over the
ideal $I'/(x)$.  Therefore we obtain $u \in \left[ \;
\overline{I'/(x) } \; \right]_1$.
\end{proof}

\smallskip

In the remainder of this section, we are going to strengthen 
Theorem \ref{HImod1} under more stringent hypotheses, 
in that we replace $\overline{E'/R'x}^{\;R'^e/R'x}$
by $\overline{E'/R'x}$, the integral closure of $E'/R'x$ in a free module of the same rank. 
Recall that this integral closure does not depend on the choice of 
embedding into a free module (of arbitrary rank even) if the ring is normal.

\begin{Theorem}\label{HImod2}
Let $R$ be a Noetherian, universally catenary, normal ring that is locally
analytically unramified. Let $E=\sum_{i=1}^n R a_i$ be a torsion-free $R$--module 
having a rank $e \geq 2$. Let
$R'=R[z_1, \ldots, z_n]$, $E'=R' \otimes_R E$, and $x=\sum_{i=1}^{n} z_{i}a_i$. 
Suppose that for
every $\mathfrak{p} \in \Spec(R)$ with $\depth R_{\mathfrak{p}} \leq
2$ either $E_{\mathfrak{p}}$ has a nontrivial free direct summand or
$R_{\mathfrak{p}}$ is regular. Then the $R'$--module $E'/R'x$ is torsion-free
and $\overline{E'/R'x}=\overline{E'}/R'x$.

\end{Theorem}

\begin{proof} By Lemma~\ref{emb0}, there exists an embedding $E
\hookrightarrow R^e$ such that $\nu((R^e/E)_{\mathfrak{p}}) \leq
e-1$ whenever $\depth R_{\mathfrak{p}} \leq 1$. By Lemma~\ref{emb1},
$R'^e/R'x$ is a  torsion-free module of rank $e-1$. Therefore 
$E'/R'x$ is torsion-free and there
exists an embedding $ R'^e/R'x \hookrightarrow R'^{e-1}$. Consider
the following diagram,

\[\xymatrix{
E'/R'x \ar@{^{(}->}[rr] \ar@{^{(}->}[dr] & & R'^e/R'x
\ar@{^{(}->}[rr] &&
R'^{e-1} \\
& \overline{E'}/R'x \ar@{^{(}->}[rr] \ar@{^{(}->}[ur] &&
\overline{\overline{E'}/R'x}\;. \ar@{^{(}->}[ur] }
\] Let $\mathfrak{E}$ denote $\overline{E'}/R'x$. By
Theorem~\ref{HImod1}, we have
\[ \mathfrak{E} = \overline{\mathfrak{E}} \cap
(R'^e/R'x) \,  . \] Therefore we obtain 
\[
\overline{\mathfrak{E}} / \mathfrak{E} \subset R'^{e-1}/(R'^e/R'x) \, .
\] Since the projective dimension of $R'^{e-1}/(R'^e/R'x)$ is at most $2$,
we have $\depth R'_{\mathfrak{P}} \leq 2$ for every
associated prime $\mathfrak{P} \in \mbox{\rm
Ass}_{R'}(\overline{\mathfrak{E}} / \mathfrak{E})$. Hence it suffices to show 
that $\mathfrak{E}_{\mathfrak{P}}$
is integrally closed for every prime ideal $\mathfrak{P}$ with
$\depth R'_{\mathfrak{P}} \leq 2$. Let $\mathfrak{p}=\mathfrak{P}
\cap R$.

We first notice that we may replace the generating sequence $a_1, \ldots, a_n$
of $E_{\mathfrak p}$ by any other generating sequence of the same length.
Indeed, any two such sequences are related by
an invertible matrix with entries in $R_{\mathfrak p}$ and the inverse 
of such a matrix defines a linear change of variables of the polynomial ring
$R_{\mathfrak p}[z_1, \ldots, z_n]$.

Now suppose that $\depth R_{\mathfrak{p}} \leq 1$. Then
$E_{\mathfrak{p}}$ is free and we may choose $a_1, \ldots, a_e$ to
be a basis and $a_{e+1}= \ldots = a_n=0$. The freeness also gives $E_{\mathfrak p}=\overline{E}_{\mathfrak{p}}$. If
$z_1, \ldots, z_e$ are not all in $\mathfrak{P}$, then 
\[\mathfrak{E}_{\mathfrak{P}} = \left(E'/R'x
\right)_{\mathfrak{P}} = \left( R'_{\mathfrak P}a_1 \oplus  \, \cdots \, \oplus R'_{\mathfrak P}a_e \right) /R'_{\mathfrak P}(z_1a_1+ \cdots +z_ea_e) \cong
\left(R'_{\mathfrak{P}} \right)^{e-1} \, .\]
Thus $\mathfrak{E}_{\mathfrak{P}}$ is
free in this case, hence integrally closed. Otherwise, since $e \geq 2$
and $\depth R'_{\mathfrak{P}} \leq 2$, we conclude that
$e=2$ and $(z_1,z_2)R'_{\mathfrak{P}}={\mathfrak P}R'_{\mathfrak P}$.
Therefore \[\mathfrak{E}_{\mathfrak{P}} =\left (R'_{\mathfrak P}a_1 \oplus R'_{\mathfrak P}a_2 \right )/R'_{\mathfrak P}(z_1a_1+z_2a_2) \cong
(z_1, z_2)R'_{\mathfrak{P}} \, ,\] which is the maximal ideal. Hence
again $\mathfrak{E}_{\mathfrak{P}}$ is integrally closed.

Next consider the case $\depth R_{\mathfrak{p}} =2$. We have $\depth R'_{\mathfrak{P}} \leq 2= \depth R_{\mathfrak{p}}$.
Since $R'$ is a polynomial ring over $R$, we see that $\mathfrak P = \mathfrak p R'$, 
hence $R'_{\mathfrak P}=R_{\mathfrak p}(z_1, \ldots, z_n)$.
Now suppose $E_{\mathfrak{p}}$ has a nontrivial free summand. In this case we
may assume that 
$a_1=(1,0)$ in $R_{\mathfrak{p}} \oplus L=E_{\mathfrak{p}}$ and that
$a_2, \ldots, a_n$ generate $L$. Since $z_1$ is a unit in $R'_{\mathfrak{P}}$, we obtain
\[\mathfrak{E}_{\mathfrak{P}}=\left(R'_{\mathfrak P} \oplus R'_{\mathfrak P} \otimes_{R_{\mathfrak p}} \overline{L}
\right)/R'_{\mathfrak P} (z_1, z_2a_2+
\cdots+z_na_n) \cong R'_{\mathfrak P} \otimes_{R_{\mathfrak p}} \overline{L} \,  . \] Therefore
$\mathfrak{E}_{\mathfrak{P}}$ is integrally closed.
Next suppose $R_{\mathfrak{p}}$ is regular, necessarily of dimension $2$.
In this case $R'_{\mathfrak P}=R_{\mathfrak p}(z_1, \ldots, z_n)$ is a 2-dimensional
regular local ring as well.
Since $R'_{\mathfrak P}=R_{\mathfrak p}(z_1, \ldots, z_n)$, it follows that
a generic Bourbaki ideal of ${\mathfrak E}_{\mathfrak P}$ with respect to  
$\left (E'/R'x \right)_{\mathfrak{P}}$ is also a generic Bourbaki ideal of $\overline{E}_{\mathfrak{p}}$
with respect to  $E_{\mathfrak p}$. Applying
Theorem~\ref{2dimRLR} twice, now shows that ${\mathfrak E}_{\mathfrak P}$ is integrally closed, as asserted. 
\end{proof}

\vspace{.0000000001in}

\begin{Corollary}\label{HIBI}
Let $R$ be a Noetherian, universally catenary, normal local ring that is
analytically unramified. Let $E$ be a finitely generated torsion-free
$R$--module with rank $e >0$. Suppose that for every
$\mathfrak{p} \in \Spec(R)$ with $\depth R_{\mathfrak{p}} \leq 2$
either $E_{\mathfrak{p}}$ has a free direct summand of rank $e-1$ or
$R_{\mathfrak{p}}$ is regular. Then a generic Bourbaki ideal of
$\overline{E}$ with respect to $E$ is integrally closed.
\end{Corollary}

\begin{proof} We may assume that $e\geq 2$. Let $E=\sum_{i=1}^n Ra_i$, 
$R''=R(z_1, \ldots, z_n)$, $E''=R''
\otimes_R E$, and $x=\sum_{i=1}^n z_ia_i$. By Theorem~\ref{HImod2},
$E''/R''x$ is torsion-free. Let $\mathfrak{P} \in \Spec(R'')$ such that
$\depth R''_{\mathfrak{P}} \leq 2$. Once we show that either
$(E''/R''x)_{\mathfrak{P}}$ has a free summand of rank $e-2$ or
$R''_{\mathfrak{P}}$ is regular, then the assertion of the corollary follows
from Theorem~\ref{HImod2} and induction on $e$. Let $\mathfrak{p} =
\mathfrak{P} \cap R$. By assumption either $E_{\mathfrak{p}}$ has a
free summand of rank $e-1$ or $R_{\mathfrak{p}}$ is regular. 
It is clear that $R''_{\mathfrak{P}}$ is regular if $R_{\mathfrak{p}}$ is regular.
Moreover, $R_{\mathfrak{p}}$ is regular if $\depth R_{\mathfrak{p}} \leq 1$.
Hence we may assume that $\depth R_{\mathfrak{p}} =2$ and
$E_{\mathfrak{p}}$ has a free summand of rank $e-1$. We write
$E_{\mathfrak{p}} = R_{\mathfrak{p}} \oplus L$, where $L$ has a 
free summand of rank $e-2$. Now the argument of the previous proof shows that
$(E''/R''x)_{\mathfrak{P}} \cong
R''_{\mathfrak{P}} \otimes_{R_{\mathfrak p}} L$, which has a free summand of rank
$e-2$. 
\end{proof}

Under suitable hypotheses, the previous theorem and corollary say that
if a module is integrally closed then the specialization by a generic 
element of a reduction is integrally closed and so is
the generic Bourbaki ideal with respect to a reduction.
We finish by proving the converse in a very general setting.

\begin{Proposition}\label{HImodconv}
Let $R$ be a Noetherian normal local ring and let $E$ be a finitely generated
torsion-free $R$--module with rank $e \geq 2$. Let $U={\sum_{i=1}^n} Ra_i$ be a reduction 
of $E$. Let $R''=R(z_1, \ldots, z_n)$, 
$E''=R'' \otimes_{R} E$, and $x={ \sum_{i=1}^n} z_ia_i$. The $R''$--module
$E''/R''x$ is torsion-free, and if this module is integrally closed, then $E$ is integrally closed.
\end{Proposition}

\begin{proof} By Lemma~\ref{emb0}, there exists an embedding $E \subset R^e$
such that $\nu((R^e/E)_{\mathfrak{p}}) \leq e-1$ whenever $\depth R_{\mathfrak{p}} \leq 1$.
This embedding induces an embedding $U \subset R^e$ such that
$\nu((R^e/U)_{\mathfrak{p}}) \leq e-1$ whenever $\depth R_{\mathfrak{p}} \leq 1$ 
because $U_{\mathfrak p}=E_{\mathfrak p}$. 
It follows from Lemma~\ref{emb1} that
$R''^e/R''x$ is torsion-free of rank $e-1$. Hence there are embeddings $E''/R''x
\subset R''^e/R''x \subset R''^{e-1}$, by which we obtain $E''/R''x
\subset \overline{E''}/R''x \subset \overline{E''/R''x} \subset
R''^{e-1}$. Hence if $E''/R''x=\overline{E''/R''x}$, then $E''=\overline{E''}$
and therefore $E=\overline{E}$. 
\end{proof}

\end{document}